\documentclass[12pt]{article}

\usepackage{graphicx}%
\usepackage{algorithm}%
\usepackage{algorithmicx}%
\usepackage{algpseudocode}%
\usepackage{listings}%
\usepackage{amsmath,amssymb,amsfonts,amsthm}%

\usepackage{color}

\usepackage{enumitem}

\newtheorem{theorem}{Theorem}

\newtheorem{prop}[theorem]{Proposition}
\newtheorem{cor}[theorem]{Corollary}
\newtheorem{con}[theorem]{Conjecture}
\newtheorem{lemma}[theorem]{Lemma}

\newcommand{\ecc}{\rm ecc}

\title{Reconstructing a graph from the boundary distance matrix}
\author{Jos\'e C\'aceres\footnotemark[1] \and Ignacio M. Pelayo\footnotemark[2]}
\date{}

\begin{document}

\maketitle

\footnotetext[1]{Departamento de Matem\'aticas, Universidad de Almer\'{\i}a, 04120 Almer\'ia, Spain. \texttt{jcaceres@ual.es}}

\footnotetext[2]{Departament de Matem\`{a}tiques, Universitat Polit\`{e}cnica de Catalunya, 08028 Barcelona, Spain. \texttt{ignacio.m.pelayo@upc.edu}}

\begin{abstract}
A vertex $v$ of a connected graph $G$ is said to be a boundary vertex of $G$ if for some other 
vertex $u$ of $G$, no neighbor of $v$ is further away from $u$  than $v$. 
The boundary  $\partial(G)$ of  $G$ is the set of all of its boundary vertices.

The boundary distance matrix $\hat{D}_G$ of a graph $G=([n],E)$ is the square matrix of order $\kappa$, being $\kappa$ the order of $\partial(G)$, such that for every $i,j\in \partial(G)$,  $[\hat{D}_G]_{ij}=d_G(i,j)$.

Given a square matrix  $\hat{B}$ of order $\kappa$, we prove under  which  conditions $\hat{B}$ is the distance matrix $\hat{D}_T$ of the set of leaves of a tree $T$, which is precisely its boundary.

We show that if $G$ is either a block graph or a unicyclic graph, then $G$ is uniquely determined by the boundary distance matrix $\hat{D}_{G}$  of $G$ and we also conjecture that this statement holds for every connected graph $G$, whenever both the order $n$ and the boundary (and thus also the boundary distance matrix) of $G$ are prefixed.

Moreover, an algorithm for reconstructing a 1-block graph (resp., a unicyclic graph) from its boundary distance matrix is given, whose time complexity in the worst case is $O(n)$ (resp., $O(n^2)$).
\end{abstract}

\noindent
\textbf{Keywords:} Boundary, Distance matrix, Block graph, Unicyclic graph, Realizability.\\
\textbf{MSC Classification:} 05C12, 05C85, 68R10.
\section{Introduction}

While typically a graph is defined by its lists of vertices and edges, significant research has been dedicated to minimizing the necessary information required to uniquely determine a graph. 
For example, various approaches include reconstructing metric graphs from density functions~\cite{dww18}, road networks from a set of trajectories~\cite{aw12}, graphs utilizing shortest paths or distance oracles~\cite{kmz18}, labeled graphs from all r-neighborhoods~\cite{mr17}, or reconstructing phylogenetic trees~\cite{bc09}. 
Of particular note is the graph reconstruction conjecture~\cite{k42,u60} which states the possibility of reconstructing any graph on at least three vertices (up to isomorphism) from the multiset of all unlabeled subgraphs obtained through the removal of a single vertex. 
Indeed, a search with the word ``graph reconstruction'' returns more than 3 million entries.

In this paper, our focus lies in the reconstruction of graphs from the distance matrix of their boundary vertices and the graph's order. 
We are persuaded that this process could hold true for all graphs, and we state it as a conjecture (see Conjecture~\ref{mainconj}). 
It is accordingly of particular interest to explore whether this conjecture holds for specific families of graphs. Our objective herein is to establish its validity for block graphs and also for unicyclic graphs.

There is a similar line of research in the continuous setting, known as the \emph{boundary rigidity problem} (introduced in \cite{g83,m81}), with can be stated as follows:
Given a compact Riemannian manifold $(M,g)$ with boundary $\partial M$, stablish under which assumptions on $\partial M$,
the geodesic distance $d_g\vert_{\partial M \times \partial M}$, uniquely determines $g$.
For further details on this topic, see \cite{pu05,suv16,u14}.

The concept of a graph's boundary was introduced by Chartrand, Erwin, Johns and Zhang in 2003~\cite{cejz03}. Initially conceived to identify local maxima of vertex distances, the boundary has since revealed a host of intriguing properties. 
It has been recognized as geodetic~\cite{chmpps06}, serving as a resolving set~\cite{hmps13},  as a strong resolving set~\cite{ryko14} and also as doubly resolving set (see Proposition \ref{prop.bound.srs}).  
Put simply, each vertex lies in the shortest path between two boundary vertices and, given any pair of vertices $x$ and $y$, there exists a boundary vertex $v$ such that either $x$ lies on the shortest path between $v$ and $y$, or vice versa. 
With such properties, it is unsurprising that the boundary emerges as a promising candidate for reconstructing the entire graph.

Graph distance matrices represent a fundamental tool for graph users, enabling the solution of problems such as finding the shortest path between two nodes. 
However, our focus here shifts mainly towards their realizability. 
That is, given a matrix (integer, positive and symmetric), we inquire whether there exists a corresponding graph where the matrix entries represent the distances between vertices. 
In 1965, Hakimi and Yau~\cite{hy65}  presented a straightforward additional condition that the matrix must satisfy to be realizable (see Theorem~\ref{distmatgen}). 
Building upon this, in 1974,  Buneman~\cite{b74} provided the matrix characterization for being the distance matrix of a tree once we know that the graph is $K_3$-free, and Graham and Pollack~\cite{gp71} computed the determinant of the distance matrix of a tree (see Theorem~\ref{detDtrees}). 
Additionally,  Howorka~\cite{h79} in 1979, formulated conditions for the distance matrix of a block graph (see Theorem~\ref{distmatblocks}), and Lin, Liu and Lu~\cite{lll15} provided the determinant of such matrices (see Theorem~\ref{detDblocks}). 
Incidentally, we use their result to derive the converse of the Graham and Pollack's theorem. 
Furthermore, we also give an algorithmic approach to the characterization of the distance matrix of a unicyclic graph.

As  previously mentioned, we are interested in the distance matrix of a graph's boundary, a submatrix of the distance matrix. 
We seek to determine the realizability of these matrices and, while we have achieved characterization in the case of trees and block graphs, a similar analysis for unicyclic graphs remains elusive.

Finally, we present a pair of algorithms for reconstructing trees and unicyclic graphs from the distance matrix of their boundary. 
In trees, the boundary corresponds to the leaves, while in unicyclic graphs, it contains the leaves along with the vertices of the cycle with degree two. 

The paper is organized as follows: this section is finished by introducing general terminology and notation. 
In Section~\ref{sec:conjecture}, we explore the notion of boundary and its relation with distance matrices. Section~\ref{blockgraphs} is devoted to the reconstruction of trees and 1-block graphs, completing first with the realizability of both the distance matrix and the boundary distance matrix of a tree and a block graph. 
In a similar way, Section~\ref{unic} undertakes the characterization of distance matrices for unicyclic graphs, followed by their reconstruction from the boundary distance matrix. 
Finally, the paper ends with a section on conclusions and open problems.

\subsection{Basic terminology}

All the graphs considered  are undirected, simple, finite and (unless otherwise stated) connected.
If $G=(V,E)$ is a graph of order $n$ and size $m$, it means that $|V|=n$ and $|E|=m$. 
Unless otherwise specified, $n\ge2$ and $V=[n]=\{1,\ldots,n\}$.

Let $v$ be a vertex of a graph $G$.
The \emph{open neighborhood} of $v$ is $N(v)=\{w \in V(G) :vw \in E\}$, and the \emph{closed neighborhood} of $v$ is $N[v]=N(v)\cup \{v\}$ .
The \emph{degree} of $v$ is $\deg(v)=|N(v)|$.
The minimum degree  (resp. maximum degree) of $G$ is $\delta(G)=\min\{\deg(u):u \in V(G)\}$ (resp. $\Delta(G)=\max\{\deg(u):u \in V(G)\}$).
If $\deg(v)=1$, then $v$ is said to be a  \emph{leaf} of $G$ and
the set and the number of leaves of $G$ are denoted by ${\cal L}(G)$ and $\ell(G)$, respectively.

Let $K_n$, $P_n$, $W_n$ and $C_n$ be, respectively, the complete graph, path, wheel and cycle of order $n$. Moreover, 
$K_{r,s}$ denotes the complete bipartite graph whose maximal independent sets are $\overline{K_r}$ and $\overline{K_s}$, respectively.
In particular, $K_{1,n-1}$ denotes the star with $n-1$  leaves.

Given a graph $G=(V,E)$ and a subset of vertices $W\subseteq V$,
the subgraph of $G$ induced by $W$, denoted by $G[W]$, has $W$ as vertex set and $E(G[W]) = \{vw \in E : v,w \in W\}$.
If $G[W]$ is a complete graph, then it is said to be a \emph{clique} of $G$.

Given a pair of vertices $u,v$ of a graph $G$, a $u-v$ \emph{geodesic} lies on a  $u-v$ shortest  path, i.e., a  path joining $u$ and $v$ of minimum order. 
Clearly,  all $u-v$ geodesics have the same length, and it is called the \emph{distance} between vertices $u$ and $v$ in $G$, denoted by $d_G(u,v)$, or simply by $d(u,v)$, when the context is clear. 
A set $W \subseteq V(G)$ is called \emph{geodetic} if any vertex of the graph is in a $u-v$ geodesic for some $u,v\in W$. 

The \emph{eccentricity} $\ecc({\it v})$  of a vertex $v$ is the distance to a farthest vertex from $v$. 
The \emph{radius} and \emph{diameter} of $G$ are respectively, the minimum and maximum eccentricity of its vertices and are denoted as ${\rm rad}(G)$ and ${\rm diam}(G)$. 
A vertex $u\in V(G)$  is a \emph{central} vertex of $G$ if $\ecc({\it u})={\rm rad}(G)$, and it is called  a  \emph{peripheral} vertex of $G$  if $\ecc({\it u})={\rm diam}(G)$.
The set of central (resp., peripheral) vertices of $G$ is called the \emph{center} (resp., \emph{periphery}) of $G$.

Let $S=\{w_1,w_2,\ldots,w_k\}$ be a set of vertices of a graph $G$.
The distance $d(v,S)$ between a vertex $v\in V(G)$ and  $S$, is the minimum of the distances between $v$ and the vertices of $S$, that is, $d(v,S)=\min\{d(v,w):w\in S\}$.
The \emph{metric representation} $r(v|S)$ of a vertex $v$ with respect to $S$ is defined as the $k$-vector 
$r(v|S)=(d_{G}(v,w_1),d_{G}(v,w_2)\ldots,d_{G}(v,w_k))$. 
A set of vertices $S$  is called \emph{resolving} if  for every pair of distinct  vertices  $x,y \in V(G)$, there exist a vertex $u\in S$ such that $d(x,u)\neq d(y,u)$, or equivalently, if $r(x|S)\neq r(y|S)$. 

Resolving sets were first introduced by Slater in~\cite{s75}, and since then, many other similar concepts have been defined, such as doubly resolving sets~\cite{chmppsw07} and 
strong resolving sets~\cite{st04, op07}. 
A set of vertices $S$ is called \emph{doubly resolving} if for every pair $x,y\in V(G)$, there exist $u,v\in S$ such that $d(x,u)-d(y,u)\neq d(x,v)-d(y,v)$. 
A set of vertices $S$ is called \emph{strong resolving} if for every pair $x,y \in V(G)$, either $x$ is in a $y-v$ geodesic or $y$ is in $x-v$ geodesic, for some vertex $v \in S$.
Clearly, every doubly (resp., strong) resolving set is also resolving, but the converse is far from being  true (see Figure~\ref{6174.34}).

A \emph{cut-vertex} is a vertex whose deletion disconnects the graph. 
A maximal subgraph of $G$ without cut-vertices is a block of $G$.
In a \emph{block graph}, every block  is a clique, or equivalently,  every cycle induces a complete subgraph.
A block of a block graph is called \emph{trivial} if it is $K_2$.
Let $K_h$ be a non-trivial block of a block graph $G$ such that $V(K_h)=\{x_1,\ldots,x_h\}$.  
For every $i\in[h]$, the connected component of $G-E(K_h)$  containing $x_i$ is called the branching graph of $x_i$ and its denoted by $G_{x_i}$.
A non-trivial block $K_h$ is called \emph{exterior} is for some vertex $x\in V(K_h)$,  $G_{x_i}$ is a tree, in which case, it is called the \emph{branching tree} of $x_i$ and it is denoted by $T_{x_i}$.
The tree $T_{x_i}$ is said to be \emph{trivial} if $V(T_{x_i})=\{x_i\}$.
A $1$-block graph is a graph containing one non-trivial exterior block.

A graph $G$ whose order and size are equal is called \emph{unicyclic}. 
These graphs contain a unique cycle that is denoted as $C_g$, where $g$ is the \emph{girth} of $G$.
The connected component of $G-E(C_g)$ containing a vertex $v\in V(C_g)$  is denoted as $T_v$ and it is called the \emph{branching tree} of $v$.
The tree $T_v$ is said to be \emph{trivial} if $V(T_v)=\{v\}$.
A vertex $v\in V(G)$ is a \emph{branching vertex} if  either $v \not\in V (C_g)$ and $\deg (v) \ge 3$ or $v \in V (C_g)$ and $\deg (v ) \ge 4$.

For further information on basic Graph Theory we refer the reader to \cite{clz16}.

\section{The conjecture}\label{sec:conjecture}

\subsection{The boundary of a graph}\label{boundary}

In this subsection, we introduce one of the essential components of our work: the boundary of a graph, which was first studied by Chartrand et al. in~\cite{cejz03}. 
A vertex $v$ of a graph $G$ is said to be a \emph{boundary vertex} 
of a vertex $u$ if no neighbor of $v$ is further away from $u$ than $v$, i.e., if for every vertex $w\in N(v)$, 
$d(u,w) \le d(u,v)$. 
The set of boundary vertices of a vertex $u$ is denoted by $\partial_{G}(u)$, or simply by $\partial(u)$, when the context is clear.
Given a  pair of vertices $u,v\in V(G)$ if $v\in \partial(u)$, then $v$ is also said to be \emph{maximally distant} from $u$.
A pair of  vertices $u,v\in V(G)$ are called \emph{mutually maximally distant}, or simply $MMD$, if both $v\in \partial(u)$ and $u\in \partial(v)$.

The \emph{boundary} of $G$, denoted by $\partial(G)$, is the set of all of its boundary vertices, i.e., 
$\partial(G)=\cup_{u \in V(G)}\partial(u)$.
Notice that, as was pointed out in \cite{ryko14}, the boundary of $G$ can also be defined as the set of $MMD$ vertices of $G$, i.e., 
$$
\partial(G)=\{v\in V(G): {\rm there \, exists} \,  u\in V(G) \, {\rm such \, that} \, u,v \,{\, \rm are} \,  MMD \}
$$

\begin{theorem}[\cite{hs07,s23}]
Let  $G$ be  a graph of order $n\ge2$ with $\kappa$ boundary vertices.
Then, $\kappa=2$ if and only if $G=P_n$.  
Moreover, $\kappa=3$ if and only if either

\begin{enumerate}[label=\rm \bf(\arabic*)]

\item $G$ is a subdivision of $K_{1,3}$; or

\item
$G$ can be obtained from $K_3$ by attaching exactly one path (of arbitrary length) to
each of its vertices.
\end{enumerate}
\label{b=3}
\end{theorem}

Also, graphs with a big $\kappa$ are well-known, as the next results show.

\begin{prop} \label{b=n} 
Given a graph $G$  of order $n$ with $\kappa$ boundary vertices,

\begin{enumerate}[label=\rm \bf(\arabic*)]

\item
If ${\rm rad}(G)={\rm diam}(G)$, then $\kappa= n$.

\item
If ${\rm diam}(G)=2$, then $n-1 \le \kappa \le n$.
\newline
Moreover, $\kappa=n-1$ if and only if $G$ contains a unique central vertex.
\end{enumerate}
\end{prop}
\begin{proof}
\begin{enumerate}[label=\rm \bf(\arabic*)]

\item
Suppose that ${\rm rad}(G)={\rm diam}(G)=d$.
Take  $u\in V(G)$ and notice that $\ecc({\it u})=d$.
Let $v\in V(G)$ such that $d(u,v)=d$.
For every vertex $w\in N(v)$, $d(u,w)\le d = d(u,v)$.
Hence, $u\in \partial(u) \subseteq \partial(G)$.

\item
If ${\rm rad}(G)={\rm diam}(G)=2$, then according to the previous item, $\kappa = n$.
Suppose that ${\rm rad}(G)=1$.
Let $V(G)=U \cup W$ such that $U$ is the set of central vertices and $W$ is the set of peripheral vertices of $G$.
Observe that $W \subsetneq \partial(G)$ and that  if $|U|=h$, then $G[U]=K_h$.
If $h\ge 2$, then every central vertex belongs to the boundary of every other vertex of $G$.
If $h=1$ and $U=\{u\}$, then for every vertex $w \in W$, $u \not\in \partial(w)$, i.e., $\partial(G)=W$, which means that $\kappa=|W|=n-1$.
\end{enumerate}
\vspace{-.5cm}
\end{proof}

\begin{cor} \label{kncn} 
Let $G$ be a graph of order $n\ge3$ with $\kappa$ boundary vertices.
\begin{enumerate}[label=\rm \bf(\arabic*)]

\item
If $G \in \{K_n, K_{r,s}, C_n\}$ and $r,s\ge2$, then $\kappa=n$.

\item
If $G \in \{W_n, K_{1,n-1}\}$, then $\kappa=n-1$.

\end{enumerate}
\end{cor}

\vspace{.07cm}
As previously mentioned, the boundary exhibits several intriguing properties, like being geodetic~\cite{chmpps06} and a resolving set~\cite{hmps13}. 
However, for the scope of this paper, its status as a strong resolving set is particularly pertinent. 
Thus, we shall now develop into this concept with some detail.

\vspace{.07cm}
That notion were first defined by Seb\H o and Tannier~\cite{st04} in 2003, and later studied in~\cite{op07}. 
They were interested in extending isometric embeddings of subgraphs into the whole graph and, to ensure that, they defined a \emph{strong resolving set} of a graph $G$ as a subset $S\subseteq V(G)$ such that for any pair $x,y\in V(G)$ there is an element $v\in S$ such that  there exists either a $x-v$ geodesic that contains $y$ or a $y-v$ geodesic containing $x$. 

\vspace{.07cm}
What is crucial for our goals is that, as a consequence of the definition, it only suffices to know the distances from the vertices of a strong resolving set to the rest of  nodes, to uniquely determine the graph. 
This issue is explored in more detail in Subsection~\ref{dm}.

\vspace{.07cm}
It was proved in \cite{ryko14} that the boundary of a graph is always a strong resolving set.
We show next that  it is also a  doubly resolving set.

\begin{prop}
The boundary $\partial(G)$ of every graph $G$ is  both a strong resolving set and a doubly resolving set.
\label{prop.bound.srs}
\end{prop}
\begin{proof}
Let $u,v\in V(G)$ such that $d(u,v)=k$ and $\{u,v\}\cap\partial(G)=\emptyset$.
So, for some vertex $w_1\in N(v)$, $d(u,w_1)=k+1$.
If $w_1\in\partial(u)$, then we are done.
Otherwise, for some vertex $w_2\in N(w_1)$, $d(u,w_2)=k+2$.

Thus, after iterating this  procedure finitely many times, say $h$ times,  we will finally  find a vertex $w_h$ such that for every vertex $w\in N(w_h)$, $d(u,w) \le d(u,w_h)=k+h$, i.e., a vertex $w_h\in\partial(G)$ and a  $u-w_h$ geodesic containing vertex $v$.
Thus, $\partial(G)$ is a strong resolving set of $G$.

Now, consider the pair $\{w_h,u\}$ and take a vertex $z_1\in N(u)$ such that $d(w_h,z_1)=d(w_h,u)+1$. 
Reasoning in the same way as before, we conclude that there is a vertex $z_{\rho}\in \partial(G)$ such  that the pair $u,v$ is in a $w_h-z_{\rho}$ geodesic. 
Hence, $\partial(G)$ is a doubly resolving set of $G$.
\end{proof}

\vspace{.07cm}
Particularly, for trees, block graphs  and unicyclic graphs, the boundary is very straightforward to characterize. 

\vspace{.17cm}
\begin{prop}[\cite{s23}] \label{boundarytrees} 
Let $T$ be a tree.
Then, $\partial(T)={\cal L}(T)$.
\end{prop}
\begin{proof}
If $u\in {\cal L}(T)$ and $N(u)=\{v\}$, then notice that $u\in \partial(v)$, and thus $u\in \partial(T)$.

Take a vertex  $u\in V(T)$ such that $\deg(u)\ge2$.
If $\{v_1,v_2\} \subseteq N(u)$ then, for every vertex $w\in V(T)$, $d(w,u) < \max \{d(w,v_1),d(w,v_2)\}$.
Hence, $u\not\in \partial(G)$.
\end{proof}

\vspace{.17cm}
\begin{prop} \label{boundaryblock} 
Let $G$ be a block graph.
If ${\cal U}(G)$ denotes the set of vertices of the blocks of order  $k\ge 3$ of degree $k-1$, then
$\partial(G)={\cal L}(G) \cup {\cal U}(G)$.
\end{prop}
\begin{proof}
If $u\in {\cal L}(G)$ and $N(u)=\{v\}$, then notice that $u\in \partial(v)$.
If $u \in {\cal U}(G)$ and $K_k$ is the clique of $G$ such that $u\in V(K_k)$, then 
for every vertex $v\in V(K_k)-u$,  $u \in \partial(v)$.

Finally, take a vertex $u \not\in {\cal L}(G) \cup {\cal U}(G)$.
If $u\in V(K_k)$, $\{v_1,v_2\} \subseteq N(u)$, $v_1\in V(K_k)$ and $v_2\not\in V(K_k)$, then 
$d(w,u) < \max\{d(w,v_1),d(w,v_2)\}$, for every vertex $w\in V(G)$.
Hence, $u\not\in \partial(G)$.
\end{proof}

\vspace{.17cm}
\begin{prop} \label{boundaryunicyc} 
Let $G$ be a unicyclic graph of girth $g$.
If ${\cal U}(G)$ denotes the set of vertices of $C_g$ of degree 2, then
$\partial(G)={\cal L}(G) \cup {\cal U}(G)$.
\end{prop}
\begin{proof}
If $u\in {\cal L}(G)$ and $N(u)=\{v\}$, then notice that $u\in \partial(v)$.
If $u \in {\cal U}(G)$ and
$v\in V(C_g)$ is a vertex such that $d(u,v) \in \lfloor \frac{g}{2} \rfloor$, then
observe that $u \in \partial(v)$.

Finally, take a vertex $u \not\in {\cal L}(G) \cup {\cal U}(G)$.
If $u\in V(C_g)$, $N(u)\cap V(C_g)=\{v_1,v_2\}$ and $v_3 \in N(u) \cap V(T_u)$, then
 $d(w,u) < \max \{d(w,v_1),d(w,v_2),d(w,v_3)\}$, for every vertex $w\in V(G)$. 
Thus, $u\not\in \partial(G)$.
If $u \in V(T_v)$ for some vertex $v\in V(C_g)$, $\deg(u)\ge 2$ and $\{v_1,v_2\} \subseteq N(u)$, then  $d(w,u) < \max\{d(w,v_1),d(w,v_2)\}$, for every vertex $w\in V(G)$.
Hence, $u\not\in \partial(G)$.
\end{proof}

\subsection{The distance matrix of a graph}\label{dm}

At this point, the other relevant element of the work is introduced: distance matrices, 
and, along with some notations, a complete characterization of both the distance matrix of a tree and the distance matrix of the leaves of a tree are provided. 
This subsection concludes by showing  a conjecture, along with some related open problems.

A square matrix $D$  is called a \emph{dissimilarity matrix} if it is symmetric, all off-diagonal entries are (strictly) positive and the diagonal entries are zeroes.
A square matrix $D$ of order $n\ge3$  is called a 
\emph{metric dissimilarity matrix} if it satisfies, for any triplet  $i,j,k \in [n]$, the  \emph{triangle inequality}: $d_{ik} \le d_{ij} + d_{jk}$.

The  \emph{distance matrix} $D_G$ of a graph $G=(V,E)$ with $V=[n]$  is the square matrix of order $n$ such that,  for every $i,j \in [n]$,  $d_{ij}=d(i,j)$.
Certainly, this matrix is a metric dissimilarity matrix.
A metric dissimilarity  matrix $D$ is called a \emph{distance matrix}
if there is a graph $G$ such that $D_G = D$.

Let $S$ be a subset of vertices of order $k$ of a graph $G=(V,E)$, with $V=[n]$. 
It is denoted by $D_{S,V}$ the submatrix of $D_G$ of order $k\times n$ such that for every $i\in S$ and for every $j\in V$, $[D_{S,V}]_{ij}=d(i,j)$.

Similarly, the so-called \emph{$S$-distance matrix} of $G$, denoted by $D^G_{S}$, or simply by $D_{S}$, when the context is clear, is the square submatrix of $D_G$ of order $k$ such that for every $i,j\in S$, 
$[D_{S}]_{ij}=d(i,j)$.
If $S=\partial(G)$, then $D_{S}$ is also denoted by $\hat{D}_G$ and it is called the \emph{boundary distance matrix} of $G$.

The next result was  stated and proved in \cite{hy65} and constitutes a general characterization of distance matrices. We include here a (new) proof, for the sake of completeness.

\begin{theorem}[\cite{hy65}]
Let  $D$ be  an integer metric dissimilarity matrix of order $n$.
Then, $D$ is a distance matrix if and only if, for every $i,j\in [n]$, if $d_{ij}>1$, then there exists an integer $k\in [n]$ such that
\begin{equation}
d_{ik}=1 \,\, {\rm and} \,\, d_{ij} = d_{ik} + d_{kj}.
\end{equation}
\label{distmatgen}
\end{theorem}
\begin{proof}
The necessity of the above condition immediately follows from the definition of distance matrix.

To prove the sufficiency, we consider the non-negative symmetric square matrix $A$ of order $n$, such that, for every pair $i,j\in [n]$, $a_{ii}\in\{0,1\}$ being $a_{ij}=1$ if and only if $d_{ij}=1$.
Let $G=([n],E)$ be the  graph  such that its adjacency matrix is $A$.
Next, we show that the distance matrix of $G$ is precisely $D$.

If ${\rm diam}(G)=d$, then $d(i,j)=p\in \{1,\ldots,d\}$, for every pair of distinct vertices $i,j\in[n]$.
If $p=1$, then clearly $d(i,j)=1$ if and only if $d_{ij}=1$.
Take $2 \le p \le d$ and suppose that, if $1 \le r \le  p-1$,  then 

$$
d(i,j)=r \,\,{\rm if \,\, and \,\, only \,\, if}\,\, d_{ij}=r
$$

Let $i,j\in [n]$ such that $d_{ij}=p$.
According to condition {\rm (1)}, take $k\in[n]$ such that $d(i,k)=d_{ik}=1$ and $p=d_{ij} = d_{ik} + d_{kj} = 1 + d_{kj}$.
This means that $d(i,j)\le p$,  since  $d_{kj}=p-1$ and $d(i,j)\le d(i,k) + d(k,j) = 1+d(k,j)$.
Hence, $d(i,j)= p$ as otherwise, according to the inductive hypothesis (1), $d_{ij}=d(i,j)<p$, a contradiction.

Conversely, let $i,j\in [n]$ such that $d(i,j)=p$.
Let $k\in N(i)$ such that $d(i,j) = 1+d(k,j)$.
Since $A$ is a metric dissimilarity matrix, $d_{ij} \le d_{ik} + d_{kj}$.
This means that $d_{ij}\le p$,  since  $d(k,j)=p-1$ and $d_{ij}\le  1+d_{kj}$.
Hence, $d_{ij}= p$ as otherwise, according to the inductive hypothesis (1), $d(i,j)=d_{ij}<p$,  a contradiction.
\end{proof}

An integer metric dissimilarity  matrix $D$ of order $n\ge 3$ is called \emph{additive} if  every subset of indices $\{i,j,h,k\} \subseteq [n]$  satisfies the so-named
\emph{four-point condition}:

$$
\left \{ 
\begin{array}{c}
d_{ij} + d_{hk} \le \max\{d_{ih} + d_{jk}, d_{ik} + d_{jh}\}  \\ \\

d_{ih} + d_{jk} \le \max\{d_{ij} + d_{hk}, d_{ik} + d_{jh}\}  \\ \\

d_{ik} + d_{jh} \le \max\{d_{ij} + d_{hk}, d_{ih} + d_{jk}\}  \\
\end{array}
\right . 
$$

Notice that every metric dissimilarity  matrix of order $n=3$ is additive, which means that the  four-point condition can be seen as a strengthened version of the triangle inequality (see \cite{b74}).

A graph $G=([n],E)$ is said to satisfy the \emph{four-point condition} if its distance matrix $D_G$ is additive, that is, if  every 4-vertex set  $\{i,j,h,k\} \subseteq [n]$ :

$$
\left \{ 
\begin{array}{c}
d(i,j) + d(h,k) \le \max\{d(i,h) + d(j,k),d(i,k) + d(j,h)\}  \\ \\

d(i,h) + d(j,k) \le \max\{d(i,j) + d(h,k),d(i,k) + d(j,h)\}  \\ \\

d(i,k) + d(j,h) \le \max\{d(i,j) + d(h,k),d(i,h) + d(j,k)\}  \\
\end{array}
\right . 
$$

As was pointed  out in \cite{b74}, these inequalities can be characterized as follows.

\begin{prop}[\cite{b74}] \label{4point}
Let $\{i,j,h,k\}$ be a 4-vertex set of a  graph $G=([n],E)$.
Then, the following statements are equivalent.

\begin{enumerate}[label=\rm \bf(\arabic*)]

\item $\{i,j,h,k\}$ satisfies  the four-point condition.

\item Among the three sums
$d(i,j) + d(h,k)$, $d(i,h) + d(j,k)$, $d(i,k) + d(j,h)$,
the two largest ones are equal. 
\end{enumerate}
\end{prop}

\subsection{The Conjecture}\label{subsec:conjecture}

The next result was  implicitly  mentioned in some papers \cite{k20,ryko14,st04} and proved in \cite{cp23}. 
This equivalence, along with the statement shown in Proposition~\ref{prop.bound.srs}, has served as an inspiration for the main conjecture of the paper that is presented at the end of this subsection.

\begin{theorem}[\cite{cp23}]
Let $S$  be a proper subset of vertices  of a graph $G=(V,E)$.
Then, the following statements are equivalent.

\begin{enumerate}[label=\rm \bf(\arabic*)]

\item $S$ is a strong resolving set.

\item $G$ is uniquely determined by the distance matrix $D_{S,V}$.

\end{enumerate}
\label{sdim.dmatrix}
\end{theorem}

\begin{figure}[ht]
\begin{center}
\includegraphics[width=0.41\textwidth]{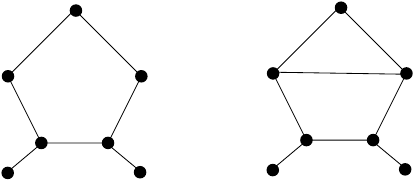}
\caption{A pair of graphs of order $7$, whose pair of leaves form a (neither doubly nor strong) resolving set.}
\label{6174.34}
\end{center}
\end{figure}

As was noticed in \cite{ryko14,st04}, this result is not true if we consider resolving sets  instead of strong resolving sets.
For example, the pair of leaves of the graphs  displayed in Figure \ref{6174.34} form, in both cases, a resolving set $S$  and also for both graphs the matrix $D_{S,V}$ is the same.

As a direct consequence of both  Theorem \ref{sdim.dmatrix} and Proposition \ref{prop.bound.srs}, the following result holds.

\begin{cor}[\cite{cp23}]
Let $G=(V,E)$ be a graph. 
Then, $G$ is uniquely determined by the distance matrix $D_{\partial(G),V}$.
\label{coro.bound.srs}
\end{cor}

\begin{figure}[ht]
\begin{center}
\includegraphics[width=0.9\textwidth]{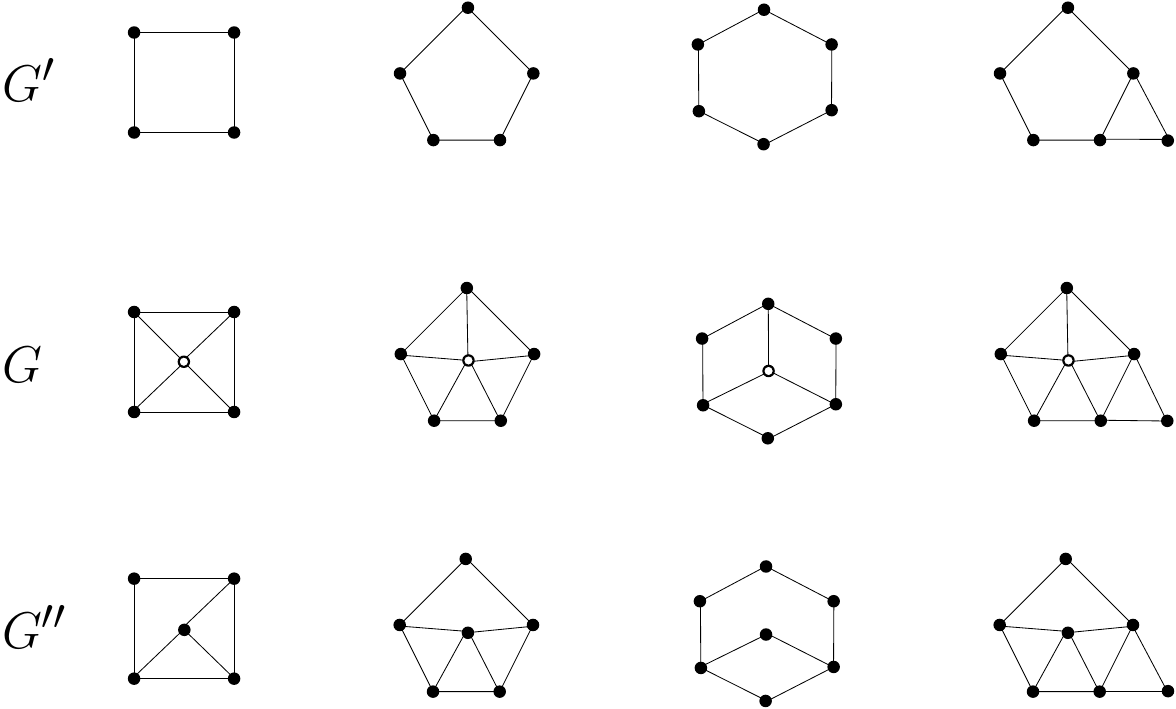}
\caption{In each column, $\hat{D'}_{G}=\hat{D}_{G}=D^{G''}_{\partial(G)}$,
$G'$ and  $G$ have the same boundary but different order, meanwhile $G$ and $G''$ have the same order but different boundary.
In all cases, the boundary is the set of black vertices.}
\label{6174.33}
\end{center}
\end{figure}

\vspace{-.1cm}
It is relatively easy to find pairs of graphs having the same boundary (that is, the same boundary distance matrix) but different order (see  Figure \ref{6174.33}, for some examples).
Having in mind all of these results and particularly  the one stated in Corollary \ref{coro.bound.srs}, we present the following conjecture.

\begin{con}
Let $\hat{B}$ an integer metric dissimilarity matrix of order $\kappa$.
Let $G=([n],E)$ be a graph such that $\hat{D}_{G}=\hat{B}$. 
If $G'=([n],E')$ is a graph such that $\hat{D}_{G'}=\hat{B}$, then $G$ and $G'$ are isomorphic.
\label{mainconj}
\end{con}

Equivalently, this conjecture can be restated as follows:

\vspace{.3cm}
\noindent
{\bf Conjecture 12.}
{\it Let $\kappa,n$ be integers such that $2\le \kappa \le n$.
Let $A, \ \hat{B}$ be integer square matrices of order $n$ and $\kappa$, respectively.
Then, there is, at most, one graph $G$ such that
$V(G)=[n]$, $D_{G}=A$, $\partial(G)=[\kappa]$ and $\hat{D}_{G}=\hat{B}$.
}

\vspace{.4cm}
Let $\kappa,n$ be integers such that $2\le \kappa \le n$.
Let $A, \ \hat{B}$ be integer square matrices of order $n$ and $\kappa$, respectively.
Let $G$ be a graph  such that $V(G)=[n]$,  $\partial(G)=[\kappa]$ and $D_{G}=A$.
We define the following graph families, denoted by ${\cal H}(\kappa)$,  ${\cal H}(n)$, and ${\cal H}(\kappa,n)$,, respectively.

\begin{itemize}

\item 
$G \in {\cal H}(\kappa)$ if it is the unique graph (up to isomorphism) such that 
$\partial(G)=[\kappa]$ and $\hat{D}_{G}=\hat{B}$

\item 
$G \in {\cal H}(n)$ if it is the unique graph (up to isomorphism) of order $n$ such that
$V(G)=[n]$ and  $D_{[\kappa]}=\hat{B}$. 

\item 
$G \in {\cal H}(\kappa,n)$ if it is the unique graph (up to isomorphism) such that
$V(G)=[n]$, $\partial(G)=[\kappa]$  and  $D_{[\kappa]}=\hat{B}$.

\end{itemize}

Notice that Conjecture \ref{mainconj} can be restated as follows.

\vspace{.3cm}
\noindent
{\bf Conjecture 12.}
{\it Every graph belongs to  ${\cal H}(\kappa,n)$.}

\vspace{.4cm}
Although it is not difficult to find graphs not belonging neither to the graph family ${\cal H}(\kappa)$ nor to the graph family ${\cal H}(n)$ (see  Figure \ref{6174.33}, for some examples), we are persuaded that, for a wide spectrum of graph classes, it is possible to obtain the whole graph $G$ from its boundary distance matrix. 
In Sections~\ref{blockgraphs} and~\ref{unic}, we  prove not only  that both the block  and the  unicyclic families belong to ${\cal H}(\kappa,n)$, but also to ${\cal H}(\kappa) \cap {\cal H}(n)$.

\section{Block graphs}
\label{blockgraphs}

This section is divided into three subsections: in the first one, we revise the main results regarding the characterization of the distance matrices of block graphs, and  we prove the converse of the result of Graham and Pollack~\cite{gp71} in Theorem~\ref{distmattrees}. 
The next subsection is devoted to determine those matrices which can be the boundary distance matrix of a block graph. 
Finally, in the last subsection, we describe and check the validity of an algorithm  to reconstruct a 1-block graph having its boundary distance matrix as the only information.

\vspace{.3cm}
\subsection{The distance matrix of a block graph}

\begin{figure}[ht]
\begin{center}
\includegraphics[width=0.6\textwidth]{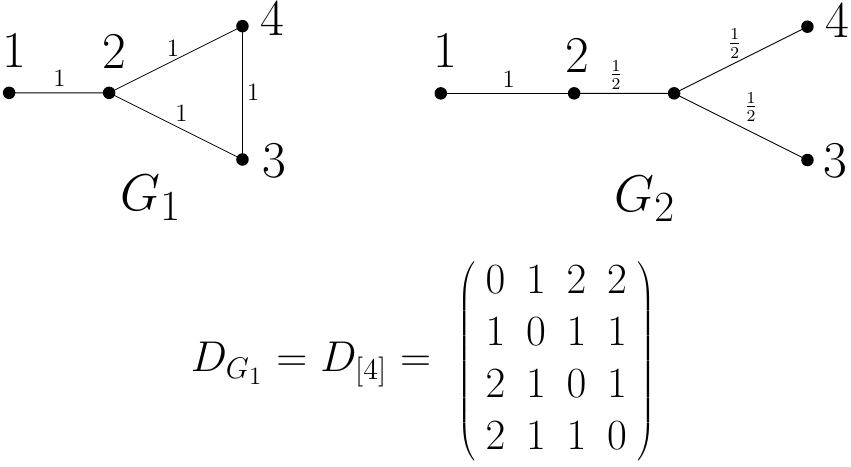}
\caption{The distance matrix $D_{G_1}$ of $G_1$  and the $[4]$-distance matrix  of  the weighted tree $G_2$ are the same.}
\label{6174.4}
\end{center}
\end{figure}

In the seminal paper \cite{b74}, Peter Buneman noticed that trees satisfy the four-point condition and also  showed that a $K_3$-free graph is a tree if and only if its distance matrix is additive.
In the same paper, it was also proved that, for every additive matrix $A$ of order $k$, there always exists a weighted tree of order $n\ge k$ containing a subset of vertices $S$ of order $k$ such that $D_S=A$ (see Figure \ref{6174.4}). A different approach based on the structure of the $4\times 4$ principal submatrices was given by Sim\~{o}es Pereira in~\cite{s69}.
In addition, it was proved in \cite{wssb77} that for every dissimilarity matrix $D$, it satisfies de four-point condition if and only if there is a unique weighted binary tree $T$ whose $\partial(T)$-distance matrix is $D$.

Starting from these results, Edward Howorka in \cite{h79} was able to characterize the family of graphs whose distance matrix is additive, i.e.,  satisfying the four-point condition.
We include next new proofs of those results for the sake of both completeness and clarity.

\begin{figure}[ht]
\begin{center}
\includegraphics[width=0.8\textwidth]{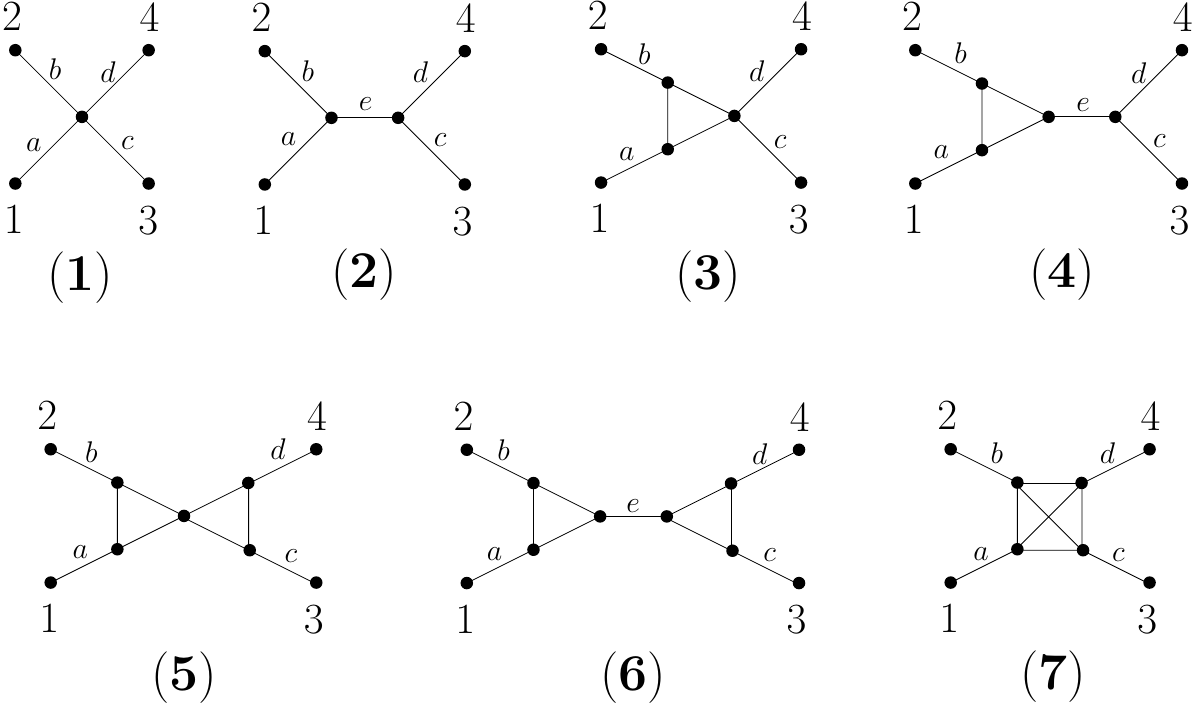}
\caption{Seven possible configurations of paths connecting 4 vertices of a block graph.}
\label{6174.5}
\end{center}
\end{figure}

\begin{prop}[\cite{h79}]
Every block graph satisfies the four-point condition.
\label{distmatblocks1}
\end{prop}
\begin{proof}
Let $S$ be a 4-vertex set of a block graph $G$, named $S=\{1,2,3,4\}$.
The only seven possible configurations of paths connecting the 4 vertices of $S$ are those shown in Figure~\ref{6174.5}.
We check that the four-condition holds in all cases.

\begin{enumerate}[label=\rm \bf(\arabic*)]

\item 
$d_{12} + d_{34}=d_{13} + d_{24}=d_{14} + d_{23}=a+b+c+d$

\item 
$d_{12} + d_{34}=a+b+c+d$

$d_{13} + d_{24}=d_{14} + d_{23}=a+b+c+d+2e$

\item 
$d_{12} + d_{34}=a+b+c+d+1$

$d_{13} + d_{24}=d_{14} + d_{23}=a+b+c+d+2$

\item 
$d_{12} + d_{34}=a+b+c+d+1$

$d_{13} + d_{24}=d_{14} + d_{23}=a+b+c+d+2e+2$

\item 
$d_{12} + d_{34}=a+b+c+d+2$

$d_{13} + d_{24}=d_{14} + d_{23}=a+b+c+d+4$

\item 
$d_{12} + d_{34}=a+b+c+d+2$

$d_{13} + d_{24}=d_{14} + d_{23}=a+b+c+d+2e+4$

\item 
$d_{12} + d_{34}=d_{13} + d_{24}=d_{14} + d_{23}=a+b+c+d+2$

\end{enumerate}
\vspace{-.5cm}
\end{proof}

The converse is proved in the next proposition.

\begin{prop}[\cite{h79}]
If $G$ satisfies the four-point condition, then it is a block graph.
\label{distmatblocks2}
\end{prop}
\begin{proof}
Let $C_h$ be an induced  cycle of  $G$ of minimum  order $h\ge4$.
Then, $h=4q+r$, with $q\ge1$ and $0\le r \le 3$.
Notice that $C_h$ is not only an induced subgraph of $G$ but also isometric.
Take a 4-vertex set $\{i,j,h,k\} \subseteq V(C_h)$ such that 
$\{d_{ij},d_{jh},d_{h,k},d_{k,i} \} \subseteq \{q,q+1\}$.
Check that $d_{ij}+d_{hk} \le 2q+2$, $d_{ik}+d_{jh} \le 2q+2$ and
$d_{ih}+d_{jk} \ge 4q$.
Hence, this 4-vertex set violate the four-point condition, which means that either $G$ is a tree or  it is a chordal graph of girth 3, i.e., the only induced cycles have length 3.

Next, suppose that $G$ is  a chordal graph of girth 3.
Take a cycle $C_p=([p],E)$ in $G$ of minimum order $p\ge4$,  such that $[p]$ is not a clique.
Notice that $p\ge 5$, since neither the cycle $C_4$ nor the diamond $K_4-e$ satisfies the four-point condition.
Let $i,j \in [p]$ such that $1 \le i <  j \le p$ and $ij\not\in E(G)$.
Notice that $d(i,j)=2$, since $C_p$ is of minimum order.
W.l.o.g. we may assume that $i=1$ and $j=3$.
Observe that for every $h \not\in \{2,p\}$ and $k \not\in \{2,4\}$,  $\{1h,3k\} \cap E(G)=\emptyset$, since $C_p$ is of minimum order.

Let $h$ be the minimum integer between $4$ and $p$  such that $2h\in E(G)$.
Then,  clearly $h=4$, since otherwise  the set $\{2,\ldots,h\}$ is an induced cycle of order at least 4, a  contradiction.
Let  $k$ be the minimum integer between $5$ and $p$  such that $2k\in E(G)$.
We distinguish  cases.

{\bf Case 1.} If $k=5$, then the subgraph induced by the set $\{2,3,4,5\}$ is the diamond $K_4-e$, a contradiction.

{\bf Case 2.} If  $p=5$ and $2k\in E(G)$, then the subgraph induced by the set $\{2,4,5,1\}$ is the cycle $C_4$, a contradiction.

{\bf Case 3.} If  $p\ge6$ and $k\ge6$, then the subgraph induced by the set $\{2,4,5,\ldots, k\}$ is the cycle $C_{k-2}$, a contradiction.

Hence, we have proved that  every cycle of $G$ induces a clique, i.e., $G$ is a block graph.
\end{proof}

Once the two implications have been proved, we can establish the theorem.

\begin{theorem}[\cite{h79}]
A  graph $G$  of order $n$ is a block graph if and only if its distance matrix $D_G$ is additive.
\label{distmatblocks}
\end{theorem}

\begin{theorem}[\cite{lll15}]
\label{detDblocks} 
Let $G$ be a block  graph on $n$ vertices and $k$ blocks $K_{n_1},\ldots,K_{n_k}$. 
Then, 
$$
\det({ D_G}) = (-1)^{n-1} \sum_{i=1}^k \frac{n_i-1}{n_i} \prod_{j=1}^k n_j
$$
\end{theorem}

In particular, as a straight consequence of the previous result,  the following theorem, proved in~\cite{gp71}, is obtained.

\begin{theorem}[\cite{gp71}]
\label{detDtrees} 
Let $T$ be a tree on $n$ vertices.Then,
$$\det({ D_T}) = (-1)^{n-1} (n-1)2^{n-2}$$
\end{theorem}

The next lemma is the crucial result that allows us to prove the characterization of distance matrices of trees by means of its determinant.

\begin{lemma}
Let  $k$ and $n$ be  integers such that $n\ge3$ and $1 \le k \le n-1$.
Let $\{ n_1,\ldots,n_k \}$ a decreasing sequence of $k$ integers  such that 
$n \ge  n_1 \ge \ldots \ge n_k \ge 2$ and $n_1+\ldots+n_k=n+k-1$.
Then,
$$
\sum_{i=1}^k \frac{n_i-1}{n_i} \prod_{j=1}^k n_j \le (n-1)2^{n-2}
$$
Moreover, the equality holds if and only if $k=n-1$ and $n_1=\ldots=n_{n-1}=2$.
\label{lemmabonito}
\end{lemma}
\begin{proof}
Let $h\in[k]$ such that $n_h\ge3$ and for every $i\in\{h+1,\ldots,k\}$, $n_i=2$.
Then,

$$
\sum_{i=1}^k \frac{n_i-1}{n_i} \prod_{j=1}^k n_j =
\Big(\frac{n_1-1}{n_1}+\ldots+\frac{n_h-1}{n_h}+\frac{k-h}{2}\Big)\cdot \prod_{j=1}^h n_j \cdot 2^{k-h}
$$

Take the $(k+1)$-sequence  $\{ n'_1,\ldots,n'_{k+1}\}=\{ n_1,\ldots, n_{h-1},n_h-1, n_{h+1},\ldots, n_k, 2 \}$. Then,

$$
\sum_{i=1}^{k+1} \frac{n'_i-1}{n'_i} \prod_{j=1}^{k+1} n'_j =
\Big(\frac{n_1-1}{n_1}+\ldots+\frac{n_{h-1}-1}{n_{h-1}}+\frac{n_h-2}{n_h-1}+\frac{k-h}{2}+\frac{1}{2}\Big)
\cdot \prod_{j=1}^{h-1} n_j \cdot (n_h-1)  \cdot  2^{k-h+1}
$$

Check that if $n_h\ge3$, then both 
$\frac{n_h-1}{n_h} < \frac{n_h-2}{n_h-1} + \frac{1}{2}$ and $n_h < (n_h-1)\cdot 2$.

Hence, 
$\sum_{i=1}^k \frac{n_i-1}{n_i} \prod_{j=1}^k n_j < \sum_{i=1}^{k+1} \frac{n'_i-1}{n'_i} \prod_{j=1}^{k+1} n'_j$.

Repeating   this  procedure iteratively, starting from the sequence $\{ n'_1,\ldots,n'_{k+1}\}$, the inequality 
$\sum_{i=1}^k \frac{n_i-1}{n_i} \prod_{j=1}^k n_j \le (n-1)2^{n-2}$ is shown, since the last sequence is the $(n-1)$-sequence: $\{2,\ldots,2\}$.
\end{proof}

As a direct consequence of Theorems \ref{distmatblocks}, \ref{detDblocks},  \ref{detDtrees} and Lemma \ref{lemmabonito}, we are able to prove the converse of Theorem~\ref{detDtrees}.

\begin{theorem}
A  graph  of order $n$ is a tree $T$ if and only if its distance matrix $D_T$ is additive and 
$\det({D_T}) = (-1)^{n-1} (n-1)2^{n-2}$.
\label{distmattrees}
\end{theorem}

\subsection{The boundary distance matrix of  a block graph}

If, in the previous subsection, we have characterized the distance matrices  of block graphs,
in this one we intend to characterize the set of metric dissimilarity matrices which are the distance matrix of the boundary of these classes of graphs.

We begin by showing  that no two (non-isomorphic) trees can have the same boundary distance matrix, a fact that was firstly noticed and proved in \cite{s62}.

\begin{theorem} 
[\cite{s62}] 
\label{thmconjtrees}
Let $T$ be a tree on $n$ vertices and $\kappa$ leaves.
Then, $T$ is uniquely determined by $\hat{D}_{T}$, the ${\cal L}(T)$-distance matrix of $T$.
\end{theorem}
\begin{proof}
We proceed by induction on ${\kappa}$. 
Clearly, the claim holds true when $\kappa=2$ since the unique tree with 2 leaves of order $n$ is the path $P_n$ and $n$  is uniquely determined by the distance between its leaves.

Let $T_{\kappa}$ be a tree with $\kappa$ leaves such that ${\cal L}(T_{\kappa}) = \{ \ell_1, \ldots, \ell_{\kappa} \}$ is the set of leaves of $T_{\kappa}$.
Assume that $\hat{D}_{T_{\kappa}}=\hat{D}_T$.
Let $\hat{D}_{{\kappa}-1}$ be the submatrix of $\hat{D}_T$ obtained by deleting the last row and column of 
$\hat{D}_T$.

By the inductive hypothesis, there is a unique tree $T_{{\kappa}-1}$ with ${\kappa}-1$  leaves such that 
$\hat{D}_{T_{{\kappa}-1}}=\hat{D}_{{\kappa}-1}$.
Hence, $T_{{\kappa}-1}$ is the subtree of $T$ obtained by deleting  the path that joins the leaf $\ell_{\kappa}$ to its exterior major vertex $w_{\kappa}$.

According to Propositions \ref{prop.bound.srs} and  \ref{boundarytrees}, ${\cal L}(T_{\kappa-1}) = \{ \ell_1, \ldots, \ell_{\kappa-1} \}$ is a doubly resolving set of $T_{\kappa-1}$.
This means that, if $d(\ell_{\kappa},w_{\kappa})=a$, then  $w_{\kappa}$ is the unique vertex in $T_{\kappa-1}$ such that 
$$
r(\ell_{\kappa}|{\cal L}(T_{\kappa-1}))  = r(w_{\kappa}|{\cal L}(T_{\kappa-1}))+(a,\stackrel{\kappa-1}{\ldots},a)
$$
Thus, $T_{\kappa}$ and $T$ are isomorphic.
\end{proof}


In \cite{z65}, the metric dissimilarity matrices which are the distance matrix of the set of leaves of a tree were characterized.

\begin{theorem}[\cite{z65}] 
Let  $\hat{B}_{\kappa}$ be  an integer metric dissimilarity matrix of order $\kappa\ge3$.
Then,  $\hat{B}_{\kappa}$ is the ${\cal L}(T)$-distance matrix of a tree $T$ if and only if it is additive and,
for every  distinct $i,j,k\in [\kappa]$,

\begin{enumerate}[label=\rm \bf(\arabic*)]

\item
$\hat{b}_{ij} < \hat{b}_{ik} + \hat{b}_{jk}$.

\item
$\hat{b}_{ij} + \hat{b}_{ik} + \hat{b}_{jk}$ is even.

\end{enumerate}
\label{distmatleavestrees}
\end{theorem}

Before approaching these pair of issues for the block graph family, we  show how to algorithmically reconstruct a tree $T$ from its ${\cal L}(T)$-distance matrix.
To this end, it is enough to notice that the proof of Theorem \ref{thmconjtrees} can be turned into an algorithm which runs in the worst case in $O(\kappa n)$ times.

\begin{algorithm}
\caption{Reconstructing-Tree}
\label{algo1}
\begin{algorithmic}[1]
\Require A matrix $\hat{D}_{T}$ of a certain tree.
\Ensure A tree $T=(V,E)$.
\State Let $\kappa$ be the order of the matrix $\hat{D}_{T}$ and let $T$ be initially a set of $\kappa$ isolated vertices $\ell_1,\ldots ,\ell_{\kappa}$;
\State Join the vertices $\ell_1$ and $\ell_2$ by a path of the length determined in $\hat{D}_{T}$;
\State Label all the vertices $u\in V$ with $r(u|\{\ell_1,\ell_2\})$, i.e., the distances from $u$ to $\{\ell_1,\ell_2\}$;
\For{$k:=3$ to $\kappa$}
		\State Compute $r(\ell_k|\{\ell_1,\ldots ,\ell_{k-1}\})$ as the distances from $\ell_k$ to $\{\ell_1,\ldots ,\ell_{k-1}\}$;
		\State \label{mainstep} Locate a vertex $u$ in $T$ and a positive integer $a$ such that $r(u|\{\ell_1,\ldots ,\ell_{k-1}\})+(a,\stackrel{k-1}{\ldots} a)=r(\ell_k|\{\ell_1,\ldots ,\ell_{k-1}\})$;
		\State Add to $T$ a path of length $a$ joining $u$ and $\ell_k$;
		\State Relabel all the vertices in $T$ with their distances to $\{\ell_1,\ldots ,\ell_k\}$;
\EndFor
\State \Return $T$.
\end{algorithmic}
\end{algorithm}

\begin{cor}
The Algorithm~\ref{algo1} runs in time $O(\kappa n)$.
\end{cor}
\begin{proof}
It is straightforward to check that the step dominating the computation is~\ref{mainstep}, and that step is repeated $O(\kappa n)$ times.
\end{proof}

\begin{cor}
Every tree $T$ of order $n$ with $\kappa$ leaves belongs not only to ${\cal H}(\kappa,n)$, but also to ${\cal H}(\kappa)$ and to ${\cal H}(n)$.
\end{cor}

\begin{theorem} 
\label{thmconjblock}
Let $G$ be a block graph on $n$ vertices and $\kappa\ge3$ boundary vertices.
Then, $G$ is uniquely determined by $\hat{D}_{G}$, the boundary distance matrix of $G$.
\end{theorem}
\begin{proof}
We proceed by induction on ${\kappa}$, the number of boundary vertices of $G$. 
For $\kappa=3$, the statement clearly holds since, according to \cite{hs07}, $G$ is either a spider or a 1-block graph whose branching trees are paths, depending on whether  $d(u_1,u_2)+d(u_1,u_3)+d(u_2,u_3)$ be either even or odd.

Let $G_{\kappa}$ be a block graph of order $n$ with $\kappa\ge4$ boundary vertices such that 
$\hat{D}_{G_{\kappa}}=\hat{D}_{G}$.
We distinguish cases.

\vspace{.2cm}
{\bf Case 1.}
There are a pair of twin vertices $u_1,u_2 \in \partial (G_{\kappa})$.
Let $G_{\kappa-1}$ the  subgraph of $G_{\kappa}$ obtained after deleting vertex $u_1$.
Notice that, according to the induction hypothesis, $G_{\kappa-1}$ is also an induced subgraph of $G$, 
since $G_{\kappa-1}$ is a block graph with $\kappa-1$ boundary vertices.
Thus, $G_{\kappa}$ and $G$ must be isomorphic since $u_1$ is in both graphs a twin of $u_2$.

\vspace{.2cm}
{\bf Case 2.}
Assume that $\partial (G_{\kappa})$  has no twins.
Let $x_1$ and $x_2$ be a pair of vertices of an exterior block of $G_{\kappa}$.
Consider its branching trees $T_{x_1}$ and $T_{x_2}$. 
If $T_{x_1}$ (resp., $T_{x_2}$) is neither trivial nor a path, recursively pruning from  $T_{x_1}$ (resp., $T_{x_2}$) 
beginning always with a leaf having maximum eccentricity, in a similar way as shown in Algorithm \ref{algo2},  as many leaves as needed until obtaining a block graph with a pair of twins.
Otherwise, delete  both $T_{x_1}-{x_1}$ and  $T_{x_2}-{x_2}$, obtaining thus a graph in which both $x_1$ and $x_2$ are twins.
In either case, we conclude that $G_{\kappa}$ and $G$ must be isomorphic since both $T_{x_1}$ and  $T_{x_2}$ are not only in $G_{\kappa}$ but also in $G$.
\end{proof}

\begin{lemma}
Let  $\hat{B}_{3}$ be  an integer metric dissimilarity matrix of order $\kappa=3$.
Then,  $\hat{B}_{3}$ is the boundary distance matrix of a block graph $G$ if and only if it satisfies the following {\bf condition}:
for every  distinct $i,j,k\in [3]$,
$\hat{b}_{ij} < \hat{b}_{ik} + \hat{b}_{jk}$.
\label{distmat3block}
\end{lemma}
\begin{proof}
If for some block graph $G$, $\hat{B}_{3}=\hat{D}_G$ is the boundary distance matrix of $G$, then it is a routine exercise to check that $\hat{B}$ satisfies the {\bf condition}.

To prove the converse, let $\hat{B}_{3}$ be  an integer metric dissimilarity matrix of order $3$:

\begin{center}
${\hat{B}_{3}}=
\left(\begin{array}{rrr} 0 & a & b \cr a &0  & c \cr b & c & 0 \end{array}\right)$
\end{center}

We distinguish cases:

\vspace{.2cm}
{\bf Case 1.}
$a+b+c$ is even.

Firstly, notice that $\min\{a,b,c\}\ge2$, since otherwise if for example $\min\{a,b,c\}=a=1$, then, according to the {\bf condition}: $b < 1+c$ and $c < 1+b$, which means that $b=c$, and thus  $a+b+c=1+2c$, a contradiction.

Consider the tree $T$ of order $n=x+y+z+1$ with 3 leaves displayed in Figure \ref{Fig5} (left) and notice that if ${\hat{B}_3}={\hat{D}_T}$, then

$$ (y+z,x+z,x+y)=(a,b,c) \Leftrightarrow (x,y,z) = \Big(\frac{b+c-a}{2},\frac{a+c-b}{2},\frac{a+b-c}{2}\Big)$$

\begin{figure}[ht]
\begin{center}
\includegraphics[width=0.8\textwidth]{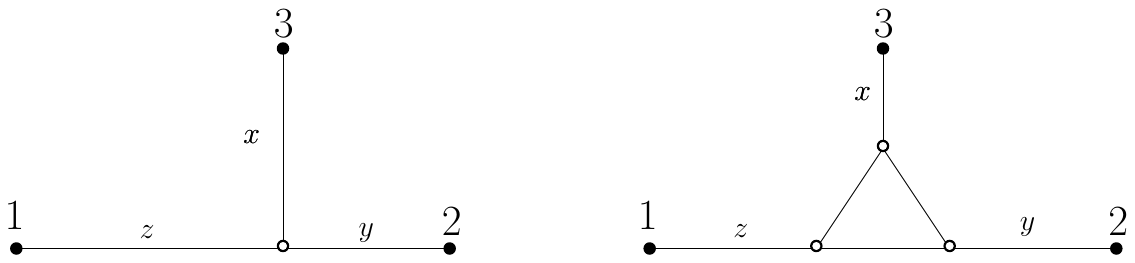}
\caption{Left: Spider of order $n=x+y+z+1$ with 3 leaves. Right: 1-block graph  of order $n=x+y+z+3$ with (at most) 3 leaves.}
\label{Fig5}
\end{center}
\end{figure}

Clearly, $x$, $y$ and $z$ are strictly positive, since $\hat{B}_3$ satisfies the {\bf condition}.
Moreover, $x$, $y$ and $z$ are integers, since  $a+b+c$ is an even integer, which means that integers $b+c-a$, $a+c-b$ and $a+b-c$ are also even.
Hence, the distance matrix of the leaves of $T$ is $\hat{B}$.

\vspace{.2cm}
{\bf Case 2.}
$a+b+c$ is odd.

Consider the 1-block graph $G$ of order $n=x+y+z+3$ with (at most) 3 leaves displayed in Figure \ref{Fig5} (right) and notice that if ${\hat{B}_3}={\hat{D}_G}$, then

{
\scriptsize
$$(y+z+1,x+z+1,x+y+1)=(a,b,c) \Leftrightarrow (x,y,z) = \Big(\frac{b+c-a-1}{2},\frac{a+c-b-1}{2},\frac{a+b-c-1}{2}\Big)$$
}

Clearly, $x$, $y$ and $z$ are  positive, since $\hat{B}_3$ satisfies the {\bf condition}.
Moreover, $x$, $y$ and $z$ are integers, since  $a+b+c$ is an odd integer, which means that integers $b+c-a$, $a+c-b$ and $a+b-c$ are also odd.
Hence, the boundary distance matrix of $G$ is $\hat{B}_3$.
\end{proof}

\begin{theorem}
\label{distmatblock}
Let  $\hat{B}_{\kappa}$ be  an integer metric dissimilarity matrix of order $\kappa\ge3$.
Then,  $\hat{B}_{\kappa}$ is the boundary distance matrix of a block graph $G$ if and only if it is additive and it satisfies the following {\bf condition}:
for every  distinct $i,j,k\in [\kappa]$,
$\hat{b}_{ij} < \hat{b}_{ik} + \hat{b}_{jk}$.
\end{theorem}
\begin{proof}
If for some block graph $G$ with $\kappa$ boundary vertices, $\hat{B}_{\kappa}=\hat{D}_T$ is the boundary distance matrix of  $G$, then it is a routine exercise to check that $\hat{B}$ is additive and satisfies the {\bf condition}.

To prove the converse, take an  integer additive matrix $\hat{B}_{\kappa}$ of order $\kappa\ge4$ satisfying the {\bf condition}.
We proceed by induction on ${\kappa}$, the order of $\hat{B}_{\kappa}$. 
Case $\kappa=3$ has been proved in Lemma \ref{distmat3block}.

Let $\hat{B}^{1}_{{\kappa}-1}$ and $\hat{B}^{\kappa}_{{\kappa}-1}$   be the matrices  obtained by deleting row (and thus also column) 1 and $\kappa$ of $\hat{B}_{\kappa}$, respectively.
Let $\hat{B}_{{\kappa}-2}$  be the matrix  obtained by deleting rows (and thus also columns) 1 and $\kappa$ of
$\hat{B}_{\kappa}$.

By the inductive hypothesis, $\hat{B}_{{\kappa}-2}$,$\hat{B}^{1}_{{\kappa}-1}$ and $\hat{B}^{\kappa}_{{\kappa}-1}$  are, respectively,  the  boundary distance matrices  of   three block graphs: $G_{\kappa-2}$, $G^{1}_{\kappa-1}$ and $G^{\kappa}_{\kappa-1}$.
Moreover, according to Theorem \ref{thmconjblock}, $G_{\kappa-2}$ is an induced subgraph of both $G^{1}_{\kappa-1}$ and $G^{\kappa}_{\kappa-1}$.

Let $G_{\kappa}$ the block graph obtained by joining $G^{1}_{\kappa-1}$ and $G^{\kappa}_{\kappa-1}$.
If ${\partial}(G^{1}_{\kappa-1}) = \{ u_2, \ldots,u_{\kappa}  \}$, 
${\partial}(G^{\kappa}_{\kappa-1}) = \{ u_1, \ldots,u_{\kappa-1}  \}$ and 
${\partial}(G_{\kappa}) = \{ u_1, \ldots, u_{\kappa-1},u_{\kappa} \}$, 
then for every $i,j\in [\kappa]$, $d(u_i,u_j)=\hat{b}_{ij}$, unless $i=1$ and $j=\kappa$.

Suppose that, for every 4-subset $\{i,j,h,k\}$, 
$\hat{b}_{ij} + \hat{b}_{hk} = \hat{b}_{ih} + \hat{b}_{jk} = \hat{b}_{ik} +\hat{b}_{jh}$.
In this case, according to Proposition \ref{distmatblocks1} (see Cases (1) and (2)),  $G_{\kappa}$ must be  either a spider with $\kappa$ legs or a 1-block graph  with  $\kappa$ branching trees, all of them being  paths (see Figure \ref{6174.5}, (1) and (7)).

Otherwise, assume w.l.o.g. that 
$\hat{b}_{12} + \hat{b}_{3\kappa} < \hat{b}_{13} + \hat{b}_{2\kappa} = \hat{b}_{1\kappa} +\hat{b}_{23}$.
Thus, $d(\ell_1,\ell_{\kappa})+d(\ell_j,\ell_{h}) = d(\ell_1,\ell_{h}) + d(\ell_j,\ell_{\kappa})$, and:
$$
d(\ell_1,\ell_{\kappa})= d(\ell_1,\ell_{h}) + d(\ell_j,\ell_{\kappa}) - d(\ell_j,\ell_{h}) 
=\hat{b}_{1h}+\hat{b}_{j\kappa}-\hat{b}_{jh}
=\hat{b}_{1\kappa},
$$
which means that $\hat{B}_{{\kappa}}$ is the distance matrix of  $G_{\kappa}$.
\end{proof}

\subsection{Reconstructing  a 1-block graph from the boundary distance matrix}

At this point, we provide a procedure to obtain a 1-block graph from its boundary distance matrix. 
The unique previous result that we need is a way to determine the leaves of the graph.

\begin{lemma}\label{lem:locateboundaryblock}
Let $G$ be a block graph. 
From  the boundary distance  matrix $\hat{D}_G$, it is possible to distinguish the vertices in ${\cal L}(G)$ from the ones in ${\cal U}(G)$.
\end{lemma}
\begin{proof}
Take a vertex  $u \in \partial(G)$.
If $u \in {\cal L}(G)$, then 
for any two distinct vertices $w_1,w_2\in \partial(G)- u$,  $d(w_1,u) + d(u,w_2)-d(w_1,w_2) \ge 2$ 
(see Figure \ref{6174.9}(1)).

If $u\in{\cal U}(G)$ and $N(u)=\{v_1,v_2$\}, consider the branching trees $T_{v_1}$ and $T_{v_2}$.
For $i=1,2$, let $w_i$ be   either a leaf of $T_{v_i}$ or the vertex $v_i$ if $T_{v_i}$, if it is trivial.
Clearly,  
$d(w_1,u) + d(u,w_2)-d(w_1,w_2) = 1$ (see Figure \ref{6174.9}(2)).
\end{proof}

Finally, Theorem~\ref{algoblock} establishes the correctness and time complexity of Algorithm~\ref{algo2}.

\vspace{.3cm}
\begin{algorithm}[ht]
\caption{Reconstructing-1Block-Recursive}\label{algo2}
\begin{algorithmic}
\Require $(\hat{B},G)$ where $\hat{B}$ is a boundary distance matrix and $G$ is a graph.
\If {$\hat{B}$ corresponds with the distance matrix of a complete graph $K_m$}
		\State \Return $(\hat{B},K_m)$
\Else
		\State Use Lemma~\ref{lem:locateboundaryblock} to distinguish the leaves in $\hat{B}$;
		\State Let $v$ be the leaf with greatest eccentricity;
		\If {$v$ has no siblings}
				\State Let $\hat{B_1}$ be the matrix $\hat{B}$ in which the row and column that correspond to the vertex $v$ have been deleted and a row and a column are added corresponding with the parent of $v$;
				\State $(\hat{B}_2,G_1)$=Reconstructing-1Block-Recursive $(\hat{B}_1,G)$;
				\State Add $v$ to $G_1$;
				\State \Return $(\hat{B},G_1)$;
		\Else { $v$ has siblings $v_1,\ldots v_k$ being $v=v_0$}
				\State Let $\hat{B_1}$ the matrix $\hat{B}$ in which the row and column that correspond to the vertex $v$ and its siblings have been deleted and a row and a column are added corresponding with the parent of $v$;
				\State $(\hat{B}_2,G_1)$=Reconstructing-1Block-Recursive $(\hat{B_1},G)$;
				\State Add $v_0,\ldots ,v_k$ to $G_1$;
				\State \Return $(\hat{B},G_1)$;		
		\EndIf
\EndIf
\end{algorithmic}
\end{algorithm}

\begin{theorem} 
\label{algoblock} 
Beginning with the boundary distance matrix $\hat{D}_{G}$ of a 1-block graph $G$, \emph{Reconstruct-1Block-Recursive} $(\hat{D}_G,\emptyset)$ obtains an isomorphic graph to $G$ in $O(n)$.
\end{theorem}
\begin{proof}
The algorithm is recursive for simplicity and to take advantage of the recursion stack for reconstructing the graph. Two parameters are involved in this process which are the distance boundary matrix $\hat{B}$ and the graph $G$. 
The matrix $\hat{B}$ is simplified by pruning one or several leaves but ensuring that the new matrix is the distance matrix of the boundary of a new graph. 
The graph $G$ plays no role in this part of the process.

The base case of the recursion occurs when $\hat{B}$ corresponds with a complete graph which is then assigned to $G$. In the backtracking process, $\hat{B}$ recuperate its previous state, with the originally pruned leaves and $G$ is actualized by adding those leaves, until the algorithm reaches the last recursive call and then $G$ is the reconstruction we are looking for.

Clearly, the algorithm pruned at least a leaf at each step, and hence the running time in the worst case is $O(n)$.
\end{proof}

\begin{cor}
Every 1-block graph $G$ of order $n$ with $\kappa$ boundary vertices belongs not only to ${\cal H}(\kappa,n)$, but also to ${\cal H}(\kappa)$ and to ${\cal H}(n)$.
\end{cor}

\section{Unicyclic graphs}\label{unic}
In this section, we move from block graphs to unicyclic graphs, i.e., those graphs containing a unique cycle.
It is divided into two subsections: one devoted to the procedure for knowing whether a matrix is the distance matrix of a unicyclic graph or not. A similar procedure works for recognizing the distance boundary matrix of a unicyclic graph.

The second one is dedicated to the process of reconstructing a unicyclic graph from its 
$\partial(G)$-distance matrix which is very similar to the analogous algorithm for 1-blocks. Incidentally, the correctness of the algorithm proves that unicyclic graphs verify Conjecture~\ref{mainconj}.

\subsection{The distance matrix of a unicyclic graph}

\begin{figure}[ht]
\begin{center}
\includegraphics[width=0.6\textwidth]{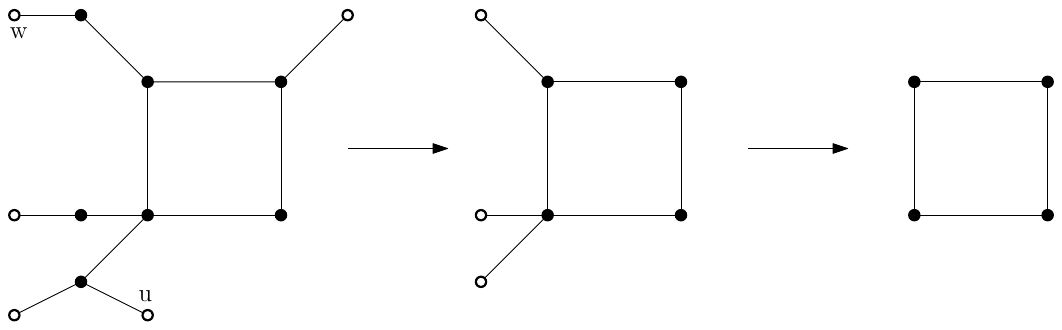}
\caption{Procedure to recognize the  distance boundary matrix of a unicyclic graph $G$.}
\label{6174.8}
\end{center}
\end{figure}

Let us focus in recognizing whether a matrix is the distance matrix of a unicyclic graph. The procedure for  checking is inductive and simple. At each step, one can delete a leaf. When there are no leaves, the resulting matrix should be one of a cycle (see Figure~\ref{6174.8}).



\begin{theorem}
\label{distmatunicyc1}
A graph $G$ is unicyclic if and only if the above procedure answers in the affirmative. 
\end{theorem}
\begin{proof}
Let $D_G$ be the distance matrix of a graph $G$. Then, a row and a column with a unique one corresponds with a leaf in the graph $G$, and if we delete that row and column, then the new matrix $D'$ is the distance matrix of the graph $G'$ obtained by deleting that leaf in $G$ (see Figure~\ref{6174.8}). 
Hence, if the final matrix of the above procedure is the distance matrix of a cycle, then it only remains to rebuild the graph to obtain a unicyclic graph.
\end{proof}

It is possible to slightly modify the previous procedure to recognize $\partial(G)$-distance matrices of unicyclic graphs. This algorithm could be recursive or iterative but in any case, we have to reduce the matrix keeping in mind that the new  matrix should be again a distance boundary matrix of a graph. In order to do that, it is only necessary to delete the leaves in a certain order. Thus, we will pick a leaf with maximum eccentricity. If that leaf has no siblings (case of the vertex $w$ in Figure~\ref{6174.8}), we delete it and substitute in the matrix for its parent which undoubtedly is a boundary vertex of the reduced graph. If the vertex is part of a bunch of siblings, then all of them are deleted and changed by its common parent (vertex $u$ in Figure~\ref{6174.8}).

It only remains a point that need to be clarify. Whereas in the distance matrix recognizing a leaf consists of determining a row or column with a unique one, in the $\partial(G)$-distance matrix, we need a different criterion for recognizing leaves which is given by the next result.

\begin{figure}[ht]
\begin{center}
\includegraphics[width=0.7\textwidth]{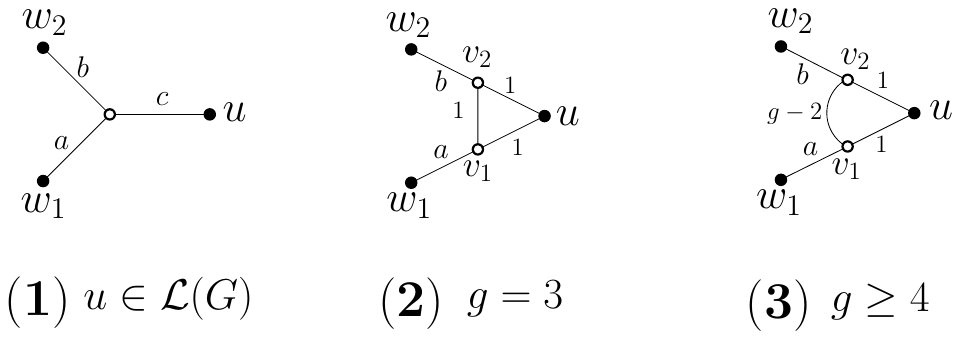}
\caption{In all cases, $u,w_1,w_2\in \partial(G)$.}
\label{6174.9}
\end{center}
\end{figure}

\begin{lemma}\label{lem:locateboundary}
Let $G$ be a unicyclic graph with $g\geq 3$. 
Given the matrix $\hat{D}_G$, it is possible to distinguish the vertices in ${\cal L}(G)$ from the ones in ${\cal U}(G)$.
\end{lemma}
\begin{proof}
Take a vertex  $u \in \partial(G)$.
If $u \in {\cal L}(G)$, then 
for any two distinct vertices $w_1,w_2\in \partial(G)- u$,  $d(w_1,u) + d(u,w_2)-d(w_1,w_2) \ge 2$ 
(see Figure \ref{6174.9}(1)).
If $u\in{\cal U}(G)$ and $N(u)=\{v_1,v_2$\}, consider the branching trees $T_{v_1}$ and $T_{v_2}$.
For $i=1,2$, let $w_i$ be   either a leaf of $T_{v_i}$ or the vertex $v_i$ if $T_{v_i}$, if it is trivial.

Clearly, if $g\ge4$ then $d(w_1,w_2) = d(w_1,u) + d(u,w_2)$ (see Figure \ref{6174.9}(3)), meanwhile that if $g=3$, then 
$d(w_1,u) + d(u,w_2)-d(w_1,w_2) = 1$ (see Figure \ref{6174.9}(2)).
\end{proof}

\begin{algorithm}[ht]
\caption{Recognizing-Unicyclic-Recursive}\label{algo1.5}
\begin{algorithmic}
\Require $\hat{D}_G$, the boundary distance matrix of a graph $G$.
\Ensure Ans=T/F depending on whether $G$ is unicyclic or not.
\State Use Lemma~\ref{lem:locateboundary} to distinguish the leaves in $\hat{D}_G$;
\If {$\hat{D}_G$ has no leaves}
		\State \Return Ans=True or False depending on $\hat{D}_G$ is the distance matrix of a cycle;
\Else
		\State Let $v$ be the leaf with greatest eccentricity;
		\If {$v$ has no siblings}
				\State Let $\hat{B}$ be the matrix $\hat{D}_G$ in which the row and column that correspond to the vertex $v$ have been deleted and a row and a column are added corresponding with the parent of $v$;
		\Else { $v$ has siblings $v_1,\ldots v_k$ being $v=v_0$}
				\State Let $\hat{B}$ the matrix $\hat{D}_G$ in which the row and column that correspond to the vertex $v$ and its siblings have been deleted and a row and a column are added corresponding with the parent of $v$;
		\EndIf
\State Ans=Recognizing-Unicyclic-Recursive $(\hat{B})$;
\EndIf
\end{algorithmic}
\end{algorithm}

\subsection{Reconstructing a unicyclic graph from the boundary distance matrix}

In this subsection, it is described the process of reconstructing a unicyclic graph $G$ from $\hat{D}_G$, the distance matrix of its boundary. 
The idea of the algorithm is the same as for 1-blocks: we prune all the leaves in a special order and the remaining graph should be a cycle graph in which we add again the leaves in reverse order. 

The algorithm implements Proposition~\ref{boundaryunicyc}, $\partial(G)={\cal L}(G)\cup {\cal U}(G)$, and Lemma~\ref{lem:locateboundary} to keep track of the leaves of the graph.



The pseudocode description is given in Algorithm~\ref{algo3}.
Finally, Theorem~\ref{algounic} establishes the correctness and time complexity of Algorithm~\ref{algo3}.

{\small
\vspace{.3cm}
\begin{algorithm}[ht]
\caption{Reconstructing-Unicyclic-Recursive}\label{algo3}
\begin{algorithmic}
\Require $(\hat{B},G)$ where $\hat{B}$ is a boundary distance matrix and $G$ is a graph.
\If {$\hat{B}$ corresponds with the distance matrix of a cycle $C_g$}
		\State \Return $(\hat{B},C_g)$
\Else
		\State Use Lemma~\ref{lem:locateboundary} to distinguish the leaves in $\hat{B}$;
		\State Let $v$ be the leaf with greatest eccentricity;
		\If {$v$ has no siblings}
				\State Let $\hat{B_1}$ be the matrix $\hat{B}$ in which the row and column that correspond to the vertex $v$ have been deleted and a row and a column are added corresponding with the parent of $v$;
				\State $(\hat{B}_2,G_1)$=Reconstructing-Unicyclic-Recursive $(\hat{B}_1,G)$;
				\State Add $v$ to $G_1$;
				\State \Return $(\hat{B},G_1)$;
		\Else { $v$ has siblings $v_1,\ldots v_k$ being $v=v_0$}
				\State Let $\hat{B_1}$ the matrix $\hat{B}$ in which the row and column that correspond to the vertex $v$ and its siblings have been deleted and a row and a column are added corresponding with the parent of $v$;
				\State $(\hat{B}_2,G_1)$=Reconstructing-Unicyclic-Recursive $(\hat{B_1},G)$;
				\State Add $v_0,\ldots ,v_k$ to $G_1$;
				\State \Return $(\hat{B},G_1)$;		
		\EndIf
\EndIf
\end{algorithmic}
\end{algorithm}
}

\begin{theorem} 
\label{algounic} 
Beginning with the boundary distance matrix $\hat{D}_{G}$ of a unicyclic graph $G$, \emph{Reconstruct-Unicyclic-Recursive} $(\hat{D}_G,\emptyset)$ obtains an isomorphic graph to $G$ in $O(n)$.
\end{theorem}
\begin{proof}
The algorithm is an evolved version of Algorithm~\ref{algo1.5} in which we added a second parameter $G$ along with the distance boundary matrix $\hat{B}$. As in the other algorithm, the matrix $\hat{B}$ is simplified by pruning one or several leaves but ensuring that the new matrix is the distance matrix of the boundary of a new graph. The graph $G$ plays no role in this part of the process.

The base case of the recursion occurs when $\hat{B}$ corresponds with a cycle graph which is then assigned to $G$. In the backtracking process, $\hat{B}$ recuperate its previous state, with the originally pruned leaves and $G$ is actualized by adding those leaves, until the algorithm reaches the last recursive call and then $G$ is the reconstruction we are looking for.

Clearly, the algorithm pruned at least a leaf at each step, and hence the running time in the worst case is $O(n)$.
\end{proof}
%
%
As a consequence, we obtain the uniqueness of the graph beginning with $\hat{D}_G$.

\begin{cor} 
\label{thmconjunicyclic}
Let $G$ be a unicyclic graph on $n$ vertices and $\kappa$ boundary vertices.
Then, $G$ is uniquely determined by $\hat{D}_{G}$, the boundary distance matrix of $G$. 
In other words, unicyclic graphs verify Conjecture~\ref{mainconj}.
\end{cor}

It is easy to check that, except for the cases with girth between 4 and 7 (see Figure \ref{6174.33}, for the cases g=4,5,6), 
every unicyclic graph $G$ of order $n$ with $\kappa$ boundary vertices belongs not only to ${\cal H}(\kappa,n)$, but also to ${\cal H}(\kappa)$ and to ${\cal H}(n)$.

\section{Conclusions and Further work}\label{cfw}

In \cite{st04}, it was firstly  implicitly  mentioned that  a resolving set $S$ of a graph $G$ is strong resolving if and only if the distance matrix $D_{S,V}$ uniquely determines the graph $G$ 
(see Theorem \ref{sdim.dmatrix}).
On the other hand, in \cite{ryko14} it was proved that the boundary $\partial(G)$ of every graph $G$ is a strong resolving set 
(see Proposition \ref{prop.bound.srs}).

Mainly having in mind this pair of results, we have  presented in Section {\bf 2}  the following conjecture.

\begin{con}
Every graph belongs to  ${\cal H}(\kappa,n)$.
\label{mainconjj}
\end{con}

In Sections 3 and 4, we have proved that if $G$ is either a block graph or a unicyclic graph, then it belongs to 
${\cal H}(\kappa,n)$, and we have also provided algorithms to recognize both 1-block and unicyclic graphs.

In addition, in Section 3, we have been able to characterize, for block graphs, both the distance matrix $D_G$ and the 
boundary distance matrix $\hat{D}_{GT}$ (see Theorems \ref{distmatblocks}, \ref{thmconjblock}  and \ref{distmatblock}).

We conclude with a list of suggested open problems. 

\begin{description}

\item [Open Problem 1:] Characterizing both  the distance matrices and the boundary distance matrices of unicyclic graphs in a similar way as it has been done for trees and for block graphs. 

\item [Open Problem 2:] Designing an algorithm for reconstructing  block graphs, in a similar way as it has been done for trees, 1-block graphs and unicyclic graphs.

\item [Open Problem 3:] Checking whether every  cactus graph belongs to ${\cal H}(\kappa)$, to ${\cal H}(\kappa)$, or at least to ${\cal H}(\kappa,n)$.

\item [Open Problem 4:] Checking whether every  split graph belongs to ${\cal H}(\kappa)$, to ${\cal H}(\kappa)$, or at least to ${\cal H}(\kappa,n)$.

\item [Open Problem 5:] Checking whether every  Ptolemaic graph belongs to ${\cal H}(\kappa)$, to ${\cal H}(\kappa)$, or at least to ${\cal H}(\kappa,n)$.

\item [Open Problem 6:] Checking whether every  graph of order $n$ with $n-1$ boundary vertices belongs to ${\cal H}(\kappa)$, to ${\cal H}(\kappa)$, or at least to ${\cal H}(\kappa,n)$.

\item [Open Problem 7:] Checking whether every   graph  of diameter 3 belongs to ${\cal H}(\kappa,n)$, and characterizing the set of graphs of diameter 3 belonging to ${\cal H}(n)$ (resp., to ${\cal H}(\kappa)$).

\end{description}




\begin{thebibliography}{99}

\bibitem{aw12}
Ahmed, M., Wenk, C. Constructing street networks from GPS trajectories. Epstein, L., Ferragina, P. (eds.) ESA 2012. LNCS 7701. Springer, Heidelberg (2012), 60--71.


\bibitem{bc09}
Brandes, U., Cornelsen, S. Phylogenetic graph models beyond trees. Discrete Appl. Math.157{\text{(10)}} (2009), 2361--2369.


\bibitem{b74}
Buneman, P. A note on the metric properties of trees. J. Combinatorial Theory Ser. B  17 (1974), 48--50.


\bibitem{chmpps06}
C\'aceres, J., Hernando, C.,  Mora, M.,  Pelayo, I. M.,  Puertas, M. L., Seara, C. On geodetic sets formed by boundary vertices. Discrete Math. 306(2) (2006), 188--198.



\bibitem{chmppsw07}
C\'aceres, J., Hernando, C.,  Mora, M.,  Pelayo, I. M.,  Puertas, M. L., Seara, C.,  Wood, D. R. On the metric dimension of Cartesian products of graphs SIAM J. Discrete Math. 21(2) (2007), 423--441.




\bibitem{cp23}
C\'aceres, J., Pelayo, I. M.. Metric Locations in Pseudotrees: A survey and new results (2023) submitted (arXiv:2307.13403v2).



\bibitem{cejz03}
Chartrand, G., Erwin, D.,  Johns, G. L.,  Zhang, P. Boundary vertices in graphs. Discrete Math.263(1-3) (2003), 25--34.



\bibitem{clz16}
Chartrand, G.,  Lesniak, L., Zhang, P. Graphs and digraphs. CRC Press, Boca Raton, FL (2016).


\bibitem{dww18}
Day, T.K., Wang, J., Wang, Y. Graph reconstruction by discrete Moore theory. In: Proceedings of the 34th International Symposium on Computational Geometry. Leibniz International Proceedings in Informatics (LIPIcs) 31 (2018), 25--34.



\bibitem{gp71}
Graham, R. L.,  Pollack, H. O. On the addressing problem for loop switching. Bell Syst. Tech. J.  50 (1971), 2495--2519.


\bibitem{g83}
Gromov, M. Filling Riemannian manifolds. J. Differential Geom. 18(1) (1983), 1--147.

\bibitem{hy65}
Hakimi, S. L., Yau, S. S. Distance matrix of a graph and its realizability. Quart. Appl. Math. 22 (1965), 305--317.



\bibitem{hs07}
Hasegawa, Y., Saito, A. Graphs with small boundary. Discrete Math. 307(14) (2007), 1801--1807.



\bibitem{hmps13}
Hernando, C.,  Mora, M.,  Pelayo, I. M., Seara, C. Some structural, metric and convex properties of the boundary of a graph. Ars Combin. 109 (2013), 267--283.



\bibitem{h79}
Howorka, E. On metric properties of certain clique graphs. J. Combin. Theory Ser. B 27(1), (1979), 67--74.




\bibitem{kmz18}
Kannan, S., Mathieu, C., Zhou, H. Graph reconstruction and verification. ACM Trans. Algorithms 14(4)(2018), 1--30.


\bibitem{k42}
Kelly, P.J. On Isometric Transformations Ph.D. thesis, University of Wisconsin (1942).




\bibitem{k20}
Kuziak, D. The strong resolving graph and the strong metric dimension of cactus graphs. Mathematics 8, 1266 (2020).





\bibitem{lll15}
Lin, H., Liu, R., Lu, X. The inertia and energy of the distance matrix of a connected graph. Linear Algebra Appl. 467 (2015), 29--39.
\bibitem{m81}
Michel, R. Sur la rigidit\'e impos\'ee par la longueur des g\'eod\'esiques. Invent Math. 65 (1981), 71-83.


\bibitem{mr17}
Mossel, E., Ross, N. Shotgun assembly of labeled graphs. IEEE Trans. Netw. Sci. Eng. 6(2) (2017), 145--157.




\bibitem{op07}
Oellermann, O. R., Peters-Fransen, J. The strong metric dimension of graphs and digraphs. Discrete Appl. Math.  155(3) (2007), 356--364.

\bibitem{pu05}
Pestov, L., Uhlmann, G. Two dimensional compact simple Riemannian manifolds are boundary distance rigid. Ann. of Math. 161(2) (2005), 1093--1110.

\bibitem{ryko14}
Rodr\'{\i}guez-Vel\'azquez, J. A., Yero, I. G., Kuziak, D., Oellermann, O. R. On the strong metric dimension of Cartesian and direct products of graphs Discrete Math. 335 (2014), 8--19.



\bibitem{st04}
S\H{e}bo, A., Tannier, E. On metric generators of graphs. Math. Oper. Res. 29(2) (2004), 383--393.



\bibitem{s69}
Sim\~{o}es Pereira, J. M. S. A note on the tree realizability of a distance matrix. J. Combinatorial Theory 6 (1969), 303--310.



\bibitem{s75}
Slater, P. J. Leaves of trees. Congr. Numer. 14  (1975), 549--559.



\bibitem{s62}
Smolenskii, Ye. A. A method for the linear recording of graphs. U.S.S.R. Comput. Math. and Math. Phys. 2(2) (1962), 396--397.


\bibitem{suv16}
Stefanov, P. Uhlmann, G., Vasy, A. Boundary rigidity with partial data. J. Amer. Math. Soc. 29(2) (2016) 299-332.

\bibitem{s23}
Steinerberger, S. The boundary of a graph and its isoperimetric inequality. Discrete Appl. Math. 338 (2023). 125--134.
\bibitem{u14}
Uhlmann, G. Inverse problems: seeing the unseen. Bull. Math. Sci. 4 (2014), 209--279.



\bibitem{u60}
Ulam, S.M. A Collection of Mathematical Problems. In:Interscience Tracts in Pure and Applied Mathematics 8. Interscience Publishers (1960).




\bibitem{wssb77}
Waterman, M. S., Smith, T. F., Singh, M., Beyer, W. A. Additive evolutionary trees. J. Theoret. Biol.64(2) (1977), 199--213.


\bibitem{z65}
Zareckiĭ, K. A. Constructing a tree on the basis of a set of distances between the hanging vertices (Russian). Uspehi Mat. Nauk 20(6) (1965), 90--92.





\end{thebibliography}
\end{document}